\documentclass{article}
\usepackage[utf8]{inputenc}
\usepackage{amsmath}
\usepackage{amssymb}
\usepackage{biblatex} %Imports biblatex package
\usepackage{caption}
\usepackage{subcaption}
\usepackage{todonotes}
\usepackage[most]{tcolorbox}
\usepackage{makecell}
\usepackage{placeins}

\usepackage[nohyperlinks]{acronym}
\usepackage{glossaries}

\newglossaryentry{alpha}{%
  name={$\alpha$},
  text={\alpha},
  description={characteristic exponent of $\alpha$-stable distribution}
}

\newglossaryentry{beta}{%
  name={$\beta$},
  text={\beta},
  description={symmetry parameter of $\alpha$-stable distribution}
}

\newglossaryentry{gamma}{%
  name={$\gamma$},
  text={\gamma},
  description={scale parameter of $\alpha$-stable distribution}
}

\newglossaryentry{delta}{%
  name={$\delta$},
  text={\delta},
  description={location parameter of $\alpha$-stable distribution}
}

\newglossaryentry{sigma}{%
  name={$\sigma$},
  text={\sigma},
  description={noise intensity}
}

\newglossaryentry{bifurcation}{%
  name={$k$},
  text={k},
  description={bifurcation parameter of Ornstein Uhlenbeck process}
}

\newglossaryentry{charf}{%
  name={$\psi$},
  text={\psi},
  description={characteristic function}
}

\usepackage{graphicx}
\graphicspath{ {} }

\title{Early warning signs of critical transitions - The $\alpha$-stable case.}
\begin{document}
%\printglossaries
%\input{aux/acro.tex}
\maketitle
Lucia S. Layritz$^1$, Ilya Pavlyukevich$^3$, Anja Rammig$^1$, Christian Kuehn$^2$

\vspace{12pt}
\begin{itemize}
    \item[$^1$] School of Life Sciences, Technical University of Munich, \\Hans-Carl-v.-Carlowitz-Platz 2, 85354 Freising, Germany
     \item[$^2$] Department of Mathematics, Technical University of Munich, \\Boltzmannstrasse 3, 85748, Garching bei München, Germany
      \item[$^3$] Institute of Mathematics, Friedrich Schiller University Jena, \\Ernst–Abbe–Platz 2, 07743 Jena, Germany
\end{itemize}

\setlength{\parindent}{0pt}

\section*{Abstract}

Statistical early warning signs can be used to identify an approaching bifurcation in stochastic dynamical systems and are now regularly employed in applications concerned with the identification of potential rapid, non-linear change or tipping points. However, the reliability of these early warning signs relies on a number of key mathematical assumptions, most notably the presence of Gaussian noise. We here show that for systems driven by non-Gaussian, $\alpha$-stable noise, the classical early warning signs of rising variance and autocorrelation are not supported by mathematical theory and their use poses the danger of spurious, false-positive results. To address this, we provide a generalized approach by introduce the scaling factor $\gamma_X$ as an alternative early warning sign. We show that in the case of the Ornstein-Uhlenbeck process, there exists a direct inverse relationship between $\gamma_{X}$ and the bifurcation parameter, telling us that $\gamma_{X}$ will increase as we approach the bifurcation. Our numerical simulations confirm theoretical results and show that our findings generalize well to non-linear, non-equilibrium systems. We thus provide a generalized, robust and applicable statistical early warning sign for systems driven by Gaussian and non-Gaussian $\alpha$-stable noise.

\section{Introduction}

Non-linear dynamical systems may exhibit rapid and irreversible state shifts upon a small change of a parameter \parencite{strogatz2015, kuehn2011}. The potential existence of such critical transitions or tipping points is a major concern in climate science and ecology \parencite{lenton2008, drijfhout2015, scheffer2001} and has been postulated for a number of climate subsystems such as the cryosphere \parencite{garbe2020, gregory2020}, hydrosphere \parencite{stocker1991, lohmann2021} or biosphere \parencite{hirota2011, rietkerk1997, chapin2005a, foley2005}.

In the case of stochastic systems, there may exist statistical early warning signs that precede the actual tipping point \parencite{scheffer2009, wiesenfeld1985}, for example, a rise in variance or autocorrelation. A range of real-world systems exhibits such signs before critical transitions \parencite{dai2013, dai2012, carpenter2011a, dakos2008} and an increase in variance and other observables has also been observed in time series data of climate elements suggested to approach tipping points \parencite{boulton2022, boers2021, boers2021a}.
	
The use of variance or autocorrelation as early warning signs sits on a robust body of mathematical theory concerned with bifurcations in stochastic dynamical systems \parencite{strogatz2015, kuehn2011, gardiner2009}. However, one key assumption of this theory is that we are working in the small noise limit of Gaussian white noise, assumptions that may not always hold in real-world applications \parencite{boettiger2012a}. Previous studies have already pointed out situations where classical early warning signs fail for other noise types \parencite{kuehn2022a, boettner2022, dutta2018}.

There is ample evidence that the assumption of Gaussian white (that is uncorrelated) noise does not hold for many climate variables including temperature, precipitation or sea level which have been shown to be correlated in time \parencite{ellerhoff2021, franzke2020, royston2018, hasselmann1976} or exhibit heavy tails, thus violating the Gaussian assumption \parencite{franzke2020, lovejoy1986}. Since climate change is expected to lead to a higher frequency of extreme events \parencite{rahmstorf2011, field2012}, the occurrence of heavy-tailed data might additionally become more frequent in the future.

One class of probability distributions that are characterized by such heavy tails are $\alpha$-stable distributions \parencite{franzke2020, nolan2020, chechkin2008}. The exception to this rule is the Gaussian (Normal) distribution which is a special case. Other known members of the class include the Cauchy (Lorenz) distribution or the Lévy distribution.  A range of real-world systems have been found which display $\alpha$-stable properties  \parencite{vandenheuvel2018, farsad2015, brockmann2010, shlesinger1995}. Notable examples in climate science and ecology include paleoclimatological temperature reconstruction from ice core data \parencite{ditlevsen1999}, foraging behavior of various animal populations \parencite{james2011, sims2008, bartumeus2005, viswanathan1996}, tree rings \parencite{lavallee2004} or the distribution of rainfall and other meteorological variables, \parencite{lovejoy1986, lovejoy1985}.

One important characteristic of non-Gaussian $\alpha$-stable distributions is that their variance and higher order moments diverge \parencite{nolan2020}. This has spawned much discussion about the applicability of $\alpha$-stable models to real data, as empirical moments will of course always be finite \parencite{lovejoy1986, lovejoy1985}. However, as \textcite{lovejoy1986} among others have pointed out, divergence simply means that we cannot expect moments to converge to a finite value but must rather assume them to continue increasing with sample size. This of course heavily challenges the use of a rising variance as an early warning sign of a tipping point. While the use of $\alpha$-stable noise in models of climate tipping is gaining traction \parencite{lucarini2022, zheng2020, yang2020, serdukova2017, ditlevsen1999a}, the impact of $\alpha$-stable driving noise on the existence and properties of early warning sign in such systems has not yet been assessed.

In this paper, we discuss the applicability and limits of classical early warning signs in the $\alpha$-stable case. We revise the basic theory of stochastic dynamical systems, $\alpha$-stable processes, and early warning signs in Section 2 and discuss potential pitfalls when applying classical early warning signs to systems driven by $\alpha$-stable noise. In Section 3 we introduce an alternative early warning sign  - the scaling factor $\gamma$ - showing that it is a natural generalization of the Gaussian variance scaling to the $\alpha$-stable case.  Lastly, in Section 4 we demonstrate the applicability of our generalized approach for simple numerical models: a linear system of Ornstein-Uhlenbeck type and a non-linear system passing through a fold bifurcation.
\section{Theoretical Background} \label{sec:theory}

\subsection{Stochastic dynamical systems}

Viewing the climate and its sub-components as a stochastic dynamical system dates back to seminal works by Hasselmann \parencite{hasselmann1976} and others, that separated the slow dynamics of climate from the fast fluctuations of weather, represented by noise. Observations are then produced by the interaction of the dynamical system with the driving noise.

We can formulate this view as a one-dimensional stochastic differential model

\begin{equation}
    dX(t) = -U'\big(X(t), k \big)dt + dN(t) = f(X(t)) \label{eq:randomdynamical}
\end{equation}

where $dX(t) = -U'\big(X (t), k \big)dt$ describes a deterministic dynamical system evolving in a potential $U$, $N(t)$ denotes a random perturbation,  $X(t)$ are realizations of the system at time $t$ and $k$ is a bifurcation parameter. 

The potential $U(X)$ can be chosen to represent any dynamical model suitable for the research task at hand. In this paper, we will consider two models: The (linear) Ornstein-Uhlenbeck process \eqref{eq:oup} and a non-linear, quadratic system \eqref{eq:fold}.

The Ornstein-Uhlenbeck process 
\begin{equation}
dX(t) = -\gls{bifurcation}X(t)dt + dN(t)  \label{eq:oup}  
\end{equation} 

originally described the movement of a particle subjected to the random influence of the surrounding fluid and friction \parencite{uhlenbeck1930}.  
The system has one fixed point at $x^* = 0$. It passes through a bifurcation at $k = 0$, where $x^*$  is stable for $\gls{bifurcation} > 0$ and unstable for $\gls{bifurcation} < 0$. Figure \ref{fig:trajectories}A gives the bifurcation diagram. 

The Ornstein-Uhlenbeck process is the most basic stochastic dynamical system and can be recovered from non-linear systems when linearizing around fixed points (as demonstrated in  \eqref{eq:linearization}).
 
The non-linear system

\begin{equation} \label{eq:fold}
    dX(t) = (k - X^2(t))dt  + dN(t)
\end{equation}

has a fold bifurcation, also at $k=0$. The system has two fixed points $X^{*^{\pm}} = \pm \sqrt{k}$  for $k > 0$ and none for $k < 0$. Figure~\ref{fig:trajectories}B gives the bifurcation diagram and stability of fixed points. 

The second important modeling choice to make is that of the random perturbation $N$. Usually, $N$ is assumed to be a Brownian motion, which is a Gaussian process $N = (N(t))_{t\geq0}$ with independent, stationary increments following a Normal (Gaussian) distribution: $N(s) - N(t) \sim \mathcal{N}(\mu = 0, \sigma_N)$. Here we want to focus on the case of symmetric $\alpha$-stable noise, which also has stationary, independent increments, but where the increments $N(s) - N(t) \sim \mathcal{S}(\alpha_N, \gamma_N)$ follow a $\alpha$-stable distribution\footnote{$\alpha$-stable processes are a subclass of Lévy processes \parencite{applebaum2009, sato1999}. For this reason, the name \textit{Lévy stable} process is also sometimes used \parencite{chechkin2008, chechkin2004}. Random walks following an $\alpha$-stable random variable are  called \textit{Lévy flights} \parencite{chechkin2008, shlesinger1995}.}. The wide class $\alpha$-stable distributions include the Gaussian as well as a range of heavy-tailed distributions.

\subsection{$\alpha$-stable random variables} \label{sec:astable}

We will briefly revise the most important properties of symmetric centered $\alpha$-stable random variables needed for out results. For this, we will follow the notation of \textcite{nolan2020} which describes $\alpha$-stable distributions $\mathcal{S}(\alpha, \beta, \gamma, \delta)$ with four parameters:

\begin{itemize}
    \item The characteristic exponent $\alpha~\in~[0,2]$, describing the tail behavior of $\mathcal{S}$
    \item The symmetry parameter $\beta~\in~[-1,~1]$, with $\beta = 0$ in the symmetric case.
    \item a scale parameter $\gamma \ge 0$
    \item a location parameter $\delta \in \mathbb{R}$, with $\delta = 0$ in the centered case.
\end{itemize}

\begin{figure}
    \centering
    \includegraphics[width=\textwidth]{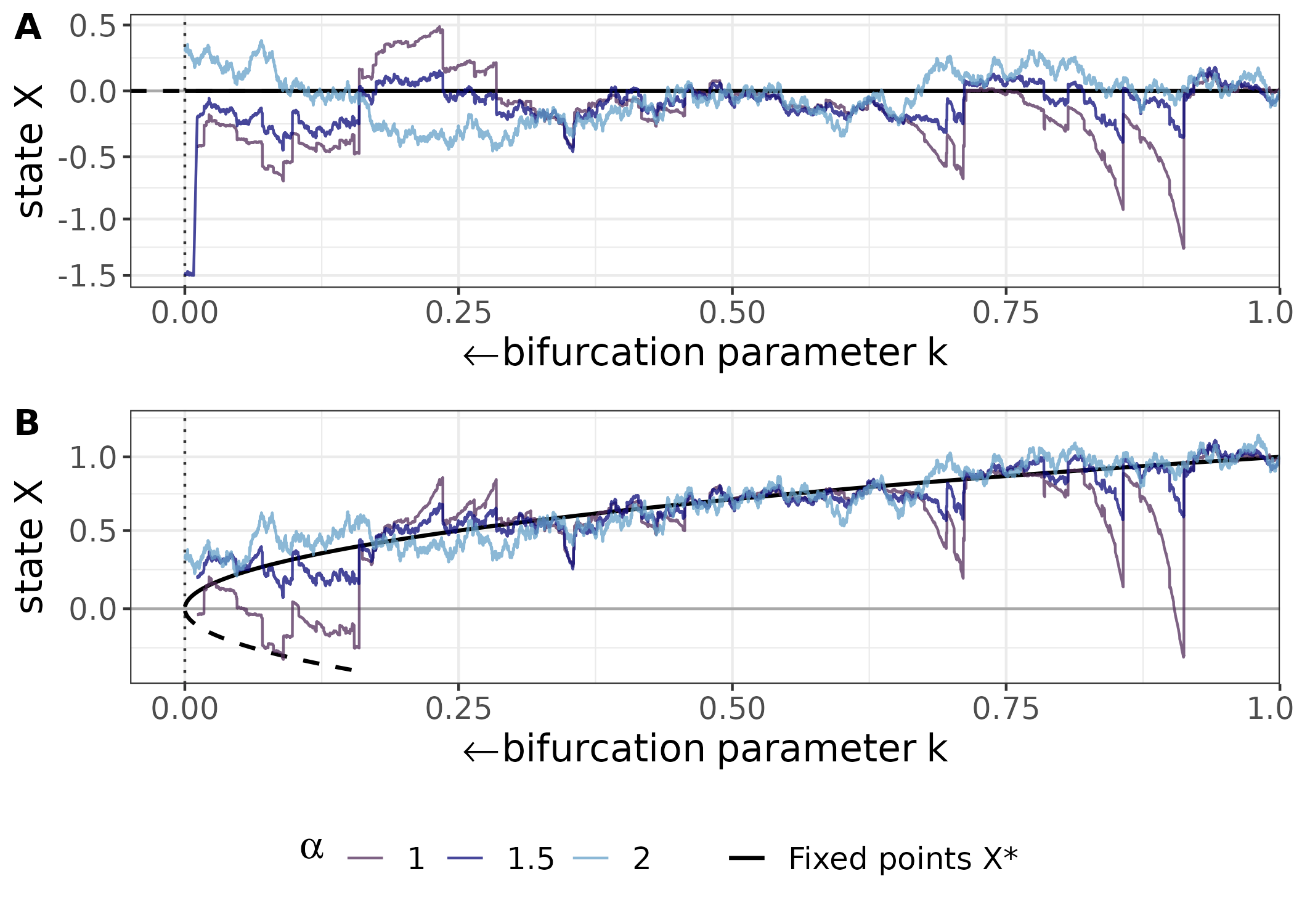}
    \caption{Bifurcation diagrams of the Ornstein-Uhlenbeck process \eqref{eq:oup} (\textbf{A}) and the fold bifurcation \eqref{eq:fold} (\textbf{B}). Stable fixed points are indicated by solid, unstable ones by dashed lines. Example trajectories for different values of $\alpha_N$ are obtained by slowly varying $k$ in the direction of the bifurcation at $k = 0$ while following the stable fixed point. All noise sequences are generated with the same random seed to ensure extreme values occur at the same time steps.}
    \label{fig:trajectories}
\end{figure}

Figure \ref{fig:theory}B illustrates the effect of the characteristic exponent $\alpha$ on the shape of the distribution and Figure~\ref{fig:trajectories}B illustrates the effect of $\alpha_N$ on trajectories of $X(t)$.

The probability density functions $f(x)$ of $\alpha$-stable random variables are in general not available. However, they can be described in terms of their characteristic function $\gls{charf}(u) = \mathbb{E}[e^{iuX}]$:

\begin{equation} \label{eq:charf_full}
    \gls{charf}(u) = e^{g(u)} \textrm{ with } g(u) = -\gls{gamma}^{\gls{alpha}}|u|^{\gls{alpha}} \textrm{ for } \beta, \delta = 0
\end{equation}

For certain special cases, however, probability density functions exist in closed form. The most important one is the Gaussian distribution, which is a special case of $\alpha$-stable distribution with  $\gls{alpha}$~=~2  and probability density function 
    \begin{equation}
        f(x) = \frac{1}{2\gls{gamma}\sqrt{\pi}} e^{{\frac{-x^2}{4\gls{gamma}^2}}}
    \end{equation}
    
The standard notation of a Gaussian density in terms of mean and variance can be recovered by substituting 
    
\begin{equation} \label{eq:gamma_variance}
  2\gls{gamma}^2 = \textnormal{Var[X]}
\end{equation}.

Other important special cases are the Cauchy distribution ($\gls{alpha}$~=~1, $\gls{beta}$~=~0) and the Lévy distribution ($\gls{alpha}$~=~$\frac{1}{2}$, $\gls{beta}$~=~1).

An important property of $\alpha$-stable distributions in the context of statistical early warning signs is that their moments $M_i = \mathbb{E}[Y^i]$ are only finite if $0 < i < \alpha$ \cite{nolan2020}. Hence the second moment and variance are not defined for all $\alpha \neq 2$ and the first moment (mean) is not defined for all $\alpha \le 1$.

\subsection{Early warning signs} \label{subsec:ews}

To construct early warning signs, we are interested in the statistical properties of $X(t)$ in relation to the bifurcation parameter $k$. We would like to reiterate that changes in these properties when approaching a bifurcation are created through the interaction of the driving noise $N$ with the dynamical system, the driving perturbation itself is assumed to remain constant. 

A range of statistical properties of $X(t)$ has been utilized as early warning signs. The most important ones, which we will focus on for the remainder of this work are variance $\textnormal{Var}[X]$ and autocorrelation $\zeta$ \parencite{wiesenfeld1985, scheffer2001, kuehn2011}, however, skewness \parencite{guttal2008} or spectral properties \parencite{bury2020} have also been proposed. 

The theory of early warning signs sits on a robust body of mathematical theory derived from the properties of the Ornstein-Uhlenbeck process \eqref{eq:oup}: For this particular system we can obtain an explicit solution (following \textcite{gardiner2009})

\begin{equation} \label{eq:oup_explicit}
    X(t) = X_0e^{-kt} +  \int_0^t e^{-\gls{bifurcation} (t-s)}dN(s)
\end{equation}

Recall that in the classical case, we assume $N$ to a Brownian motion with increments $N(s) - N(t) \sim \mathcal{N}(\mu = 0, \sigma_N)$. In this case, $X(t)$ will be normally distributed with mean $\mu = 0$ as well. We can obtain the full probability density $p(X, t)$ from Eq. \eqref{eq:oup_explicit} directly or via the Fokker-Planck-Equation

\begin{equation}
    \frac{\partial p(X,t)}{\partial t} = kX \frac{\partial p(X, t)}{\partial X} + \frac{1}{2} \frac{\partial^2p(X, t)}{\partial X^2}
\end{equation}

 to obtain the variance

\begin{equation} \label{eq:oup_var}
 \textnormal{Var}[X(t)] = (\textnormal{Var}[X_0] - \frac{\sigma_N^2}{2k})e^{-2\gls{bifurcation}t} + \frac{\sigma_N^2}{2\gls{bifurcation}} \stackrel{t \rightarrow \infty}{=} \frac{\sigma_N^2}{2\gls{bifurcation}}  
\end{equation}

\begin{figure}
    \centering
    \includegraphics[width = \textwidth]{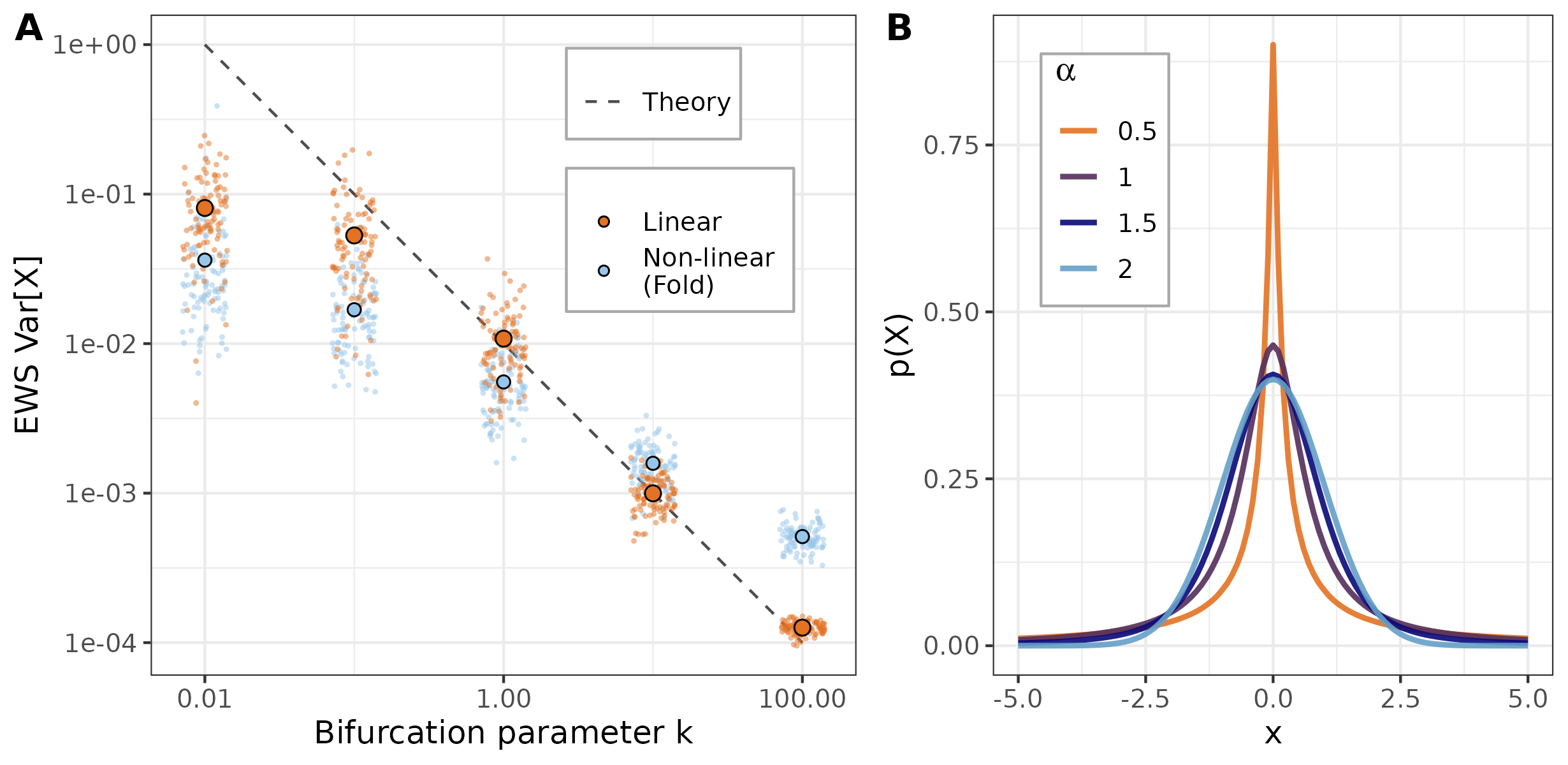}
    \caption{\textbf{A}: The classical early warning sign: Relationship between bifurcation parameter $k$ and  $\textnormal{Var}[X]$ in the Gaussian case for the Ornstein-Uhlenbeck process \eqref{eq:oup} (orange) and the fold bifurcation \eqref{eq:fold} (blue). Dashed line gives theoretical result \eqref{eq:oup_var}. \textbf{B}: Probability density functions of $\alpha$-stable random variables for different values of $\alpha$ ($\beta = 0$, $\gamma = 1$, $\delta = 0$).}
    \label{fig:theory}
\end{figure}

Assuming a deterministic initial condition ($\textnormal{Var}[X_0]$ = 0) and stationarity ($t \rightarrow \infty$), we find that the variance scales with $\frac{1}{2k}$ and hence increases as the system approaches the stable-to-unstable transition ($\gls{bifurcation}~\rightarrow~0^+$), as shown in Figure~\ref{fig:theory}A.

In non-linear systems, one would typically linearize around the steady state of interest to again obtain a linear system of the form of Eq. \eqref{eq:oup} \parencite{boers2021, boettiger2012a}. In the case of the fold bifurcation \eqref{eq:fold} we expand the right-hand side around $X^{*^{+}}$ 

\begin{equation} \label{eq:linearization}
    f(X) = f(X^{*^{+}}) + f'(X^{*^{+}})(X - X^{*^{+}}) + \mathcal{O}(|X -X^{*^{+}}|^2)
\end{equation}

After substituting  $X^{*^{+}} = \sqrt{k}$, we can rearrange to obtain a new Ornstein-Uhlenbeck process 

\begin{equation}
    dY \approx -\kappa Ydt  + dW_Y 
\end{equation}

with $Y = X - X^{*^{+}}$ and $\kappa = 2\sqrt{k}$. Hence we can expect the system to still follow relationship \eqref{eq:oup_var} when close to $X^{*^{+}}$.

The auto-correlation 

\begin{equation}
    \zeta(t_1, t_2) =  \frac{E[(X(t_1) - E[X(t_1))(X(t_2) - E[X(t_2))]}{\sqrt{\textnormal{Var}[X(t_1)]}\sqrt{(\textnormal{Var}[X(t_2)]}}
\end{equation}

follows from that, as it is a function of the first and second moment and hence mean and variance.

The result \eqref{eq:oup_explicit} also holds for an $\alpha$-stable noise process, in which case $X$ will also be $\alpha$-stable with $\alpha = \alpha_N$. 

However, as stated in Section \ref{sec:astable}, $N$ will not possess a finite variance in this case. \textcite{chechkin2004} show, that for systems of type \eqref{eq:randomdynamical} with $\alpha$-stable driving noise $N$ and $ U(X)$ of order $\frac{|x|^c}{c}$, Var$[X]$ will only be finite if $c > 4 - \alpha$ \parencite{chechkin2004}. Only then is the dynamical potential steep enough to sufficiently confine the noise. In the case of an Ornstein-Uhlenbeck process $c=2$ and therefore $c > 4 - \alpha$ does not hold with the exception of $\alpha = 2$. The same is true for the fold bifurcation $x=3$. In the case of more complex systems such as a double-well potential, the global variance may exist. Nevertheless, when we apply linearization as in equation \eqref{eq:linearization}, the local existence of variance is lost.

\begin{figure}
    \centering
    \includegraphics[width=\textwidth]{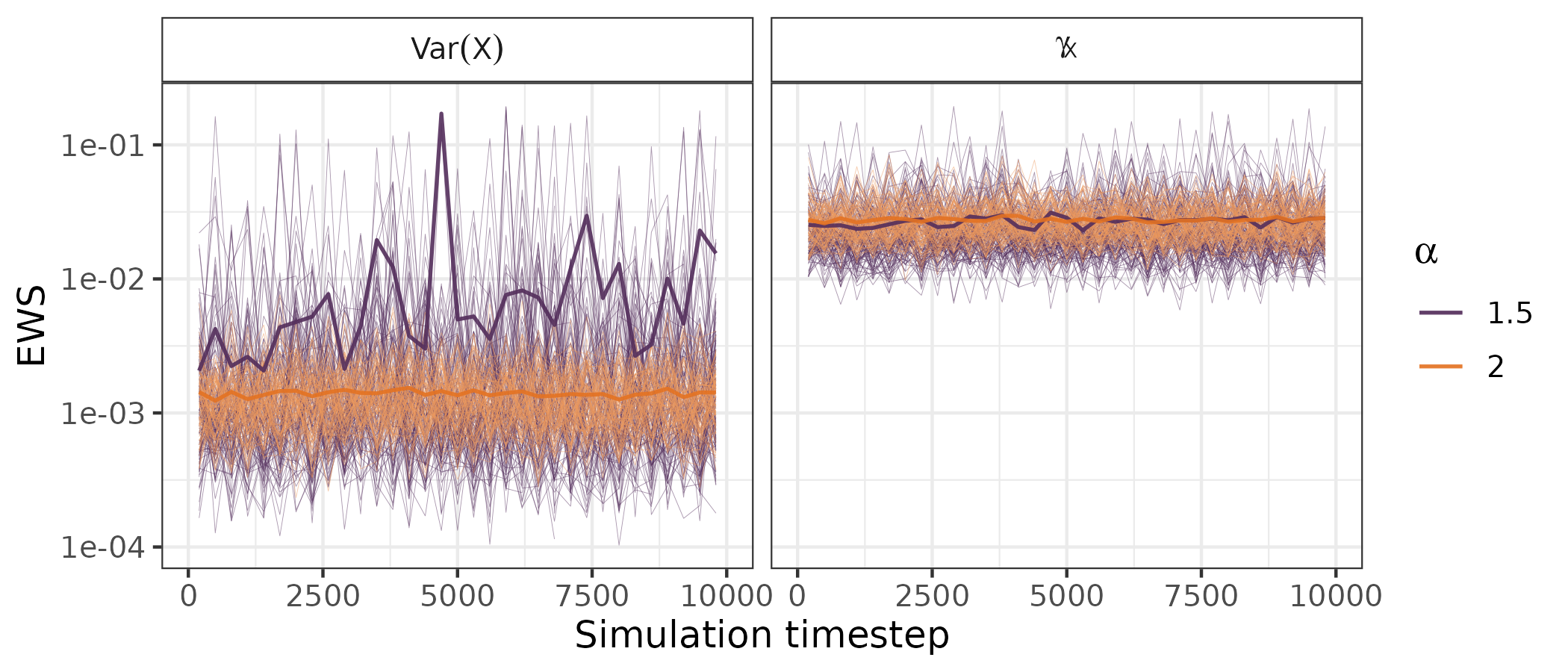}
    \caption{Illustration of non-converging variance. In the Gaussian case ($\alpha$ = 2), the variance converges to a final value once the simulation has reached equilibrium due to the central limit theorem (orange trajectories in the left panel). In the non-Gaussian, $\alpha$-stable case, the variance of $X$ does not converge to a finite value and may thus exhibit a spurious increase even for constant $k$, giving rise to a false-positive early warning sign (purple trajectories in the left panel). In contrast, $\gamma_X$ converges for all $\alpha$ (right panel). Thin lines are individually simulated trajectories, bold lines average over all 100 trajectories, $k$ = 1. Note the log scale on the y-axis. Simulation setup mirrors Figure ~\ref{fig:gamma_neq}, $\gamma_X$ is estimated every 300 simulation steps.}
    \label{fig:variance_convergance}
\end{figure}

This implies that the classical theory of early warning signs relying on linearization as laid out in this section is not valid for $\alpha$-stable systems. On the contrary, as we cannot ensure the variance to converge to a finite value, there is always the danger of misinterpreting resulting spurious increases as an early warning sign (see left panel of Figure~\ref{fig:variance_convergance} for an illustration). 

Therefore, where $\alpha$-stable systems might occur, we are in need of a different indicator that is robust again violating the Gaussian assumption.

\FloatBarrier
\newpage

\section{An early warning indicator for $\gls{alpha}$-stable systems} \label{sec:ews}

To address this caveat, we propose the scaling parameter $\gls{gamma}$ as an alternative, robust early warning sign that is applicable to Gaussian and $\alpha$-stable systems, easy to calculate in practical applications and, as we will show in the following, firmly grounded in mathematical theory. 

We start with the $\gls{alpha}$-stable process of Ornstein-Uhlenbeck $type$

\begin{equation}
dX(t) = -\gls{bifurcation}X(t)dt + dN(t)   \label{eq:stableoup}  
\end{equation}  

where $N(t)$ is a symmetric, $\alpha$-stable process with characteristic function

\begin{equation} \label{eq:charfL}
    \gls{charf}_{N}(u) = \mathbb{E}[e^{iuN(t)}] = e^{-\gls{gamma}_{N}^{\gls{alpha}}|u|^{\gls{alpha}}} 
\end{equation}

Recall the solution of the Ornstein-Uhlenbeck process \eqref{eq:oup_explicit}, which also holds in the $\alpha$-stable case

\begin{equation} \label{eq:stableoupsolution} 
    X(t) = X_0e^{-\gls{bifurcation}t} +  \int_0^t e^{-\gls{bifurcation} (t-s)}dN(s)
\end{equation}  

We know that $X(t)$ will also be $\alpha$-stable and thus can formulate its characteristic function

\begin{align} \label{eq:charfX}
    \gls{charf}_{X(t)}(u) = \mathbb{E}[e^{iuX(t)}] = e^{-\gls{gamma}_{X}^{\gls{alpha}}(t)|u|^{\gls{alpha}}} 
\end{align}

Combining Eq. \eqref{eq:stableoupsolution} into Eq. \eqref{eq:charfX} and re-arranging (see Appendix for a full derivation), we obtain

\begin{align} \label{eq:charfX_solution}
    \gls{charf}_{X(t)}(u) = e^{- \frac{\gls{gamma}_{N}^{\gls{alpha}}}{\alpha k} (1 - e^{-\alpha k t}) |u|^{\alpha}}  
\end{align}

This form allows us to retrieve the exact parameters determining the properties of $X(t)$. Comparing Eq. \eqref{eq:charfX_solution} to Eq. \eqref{eq:charfL}, we see that the random variable $X(t)$ is indeed again $\alpha$-stable with $\alpha_X = \alpha_N$ and has a scaling parameter

\begin{equation} \label{eq:gammaX}
    \gls{gamma}_{X} = \gls{gamma}_N \sqrt[\gls{alpha}]{\frac{(1 - e^{-\alpha k t})}{\gls{alpha}\gls{bifurcation}}} \stackrel{t \rightarrow \infty}{=}  \gls{gamma}_N \sqrt[\gls{alpha}]{\frac{1}{\gls{alpha}\gls{bifurcation}}}
\end{equation}

We thus find a direct relationship between $\gls{gamma}_{X}$ and the bifurcation parameter $k$, which tells us that $\gls{gamma}_{X}$ will increase as we approach the bifurcation (decreasing $k$). Based on this relationship, we are able to utilize $\gamma_X$ as an early warning sign of that bifurcation. 

This is indeed a generalization of the the variance scaling found in the Gaussian case. Recall that for Gaussian $\alpha$-stable variables $\gls{alpha} = 2$ and  $2\gls{gamma}_i^2 = \textnormal{Var}[X] = \sigma_X^2$ \eqref{eq:gamma_variance}. Substituting Eq. \eqref{eq:gammaX} into the latter gives us

\begin{equation}
   \textnormal{Var}[X] \stackrel{\eqref{eq:gamma_variance}}{=}  2\gls{gamma}^2_X \stackrel{\eqref{eq:gammaX}}{=} 2 (\gls{gamma}_N\sqrt{\frac{1}{2\gls{bifurcation}}})^2 = \frac{\gls{gamma}_N^2}{k} \stackrel{\eqref{eq:gamma_variance}}{=} \frac{\sigma^2}{2k},
\end{equation}

recovering Eq. \eqref{eq:oup_var}.

\section{Numerical Simulations} \label{sec:simulations}

We perform a range of numerical simulations to confirm our results and to illustrate the applicability of our proposed indicator $\gls{gamma}_{X}$. As \eqref{eq:gammaX} gives the solution in the long-term limit, we first perform equilibrium simulations for both systems \eqref{eq:oup}  and \eqref{eq:fold} over a range of values for $k$. In a second step, we then estimate $\gls{gamma}_{X}$ from a single trajectory while slowly moving $k$ towards the bifurcation, as one would in actual applications (non-equilibrium simulations). All simulations were performed for $\alpha$ = \{2, 1.8, 1.5, 1.3\}. We chose to focus on this range, as it is what typically occurs in real and simulated applications.

We discretize and simulate with the following Euler-Maruyama scheme \parencite{higham.2001, ditlevsen1999, samorodnitsky2017}
\begin{equation}
    X_{i} = X_{i-1} - U'(X_{i-1}) \Delta{t}+ \sqrt[\gls{alpha}]{\Delta{t}}N_{i}
\end{equation} 

where ${N_i}$ are i.i.d random variables and

\begin{equation}
    \mathbb{E}e^{iuN_k} = e^{-\gls{gamma}_{N}^{\gls{alpha}}|u|^{\gls{alpha}}}
\end{equation}

We chose $\gls{gamma}_{N} = 0.1$ and $\Delta{t} = 0.004$ and  initiated all simulations at $X_0 = 0.5$, to be in the vicinity but not at the stable state. 

Since we have more than one fixed point in the non-linear case, trajectories might escape the basin of attraction of the stable fixed point. We therefore stopped a simulation if 

\begin{equation} \label{eq:Xcritical}
    X_{i} < -\sqrt{k} - \frac{k}{10}
    \end{equation}

\begin{figure}
    \centering
    \includegraphics[width=\textwidth]{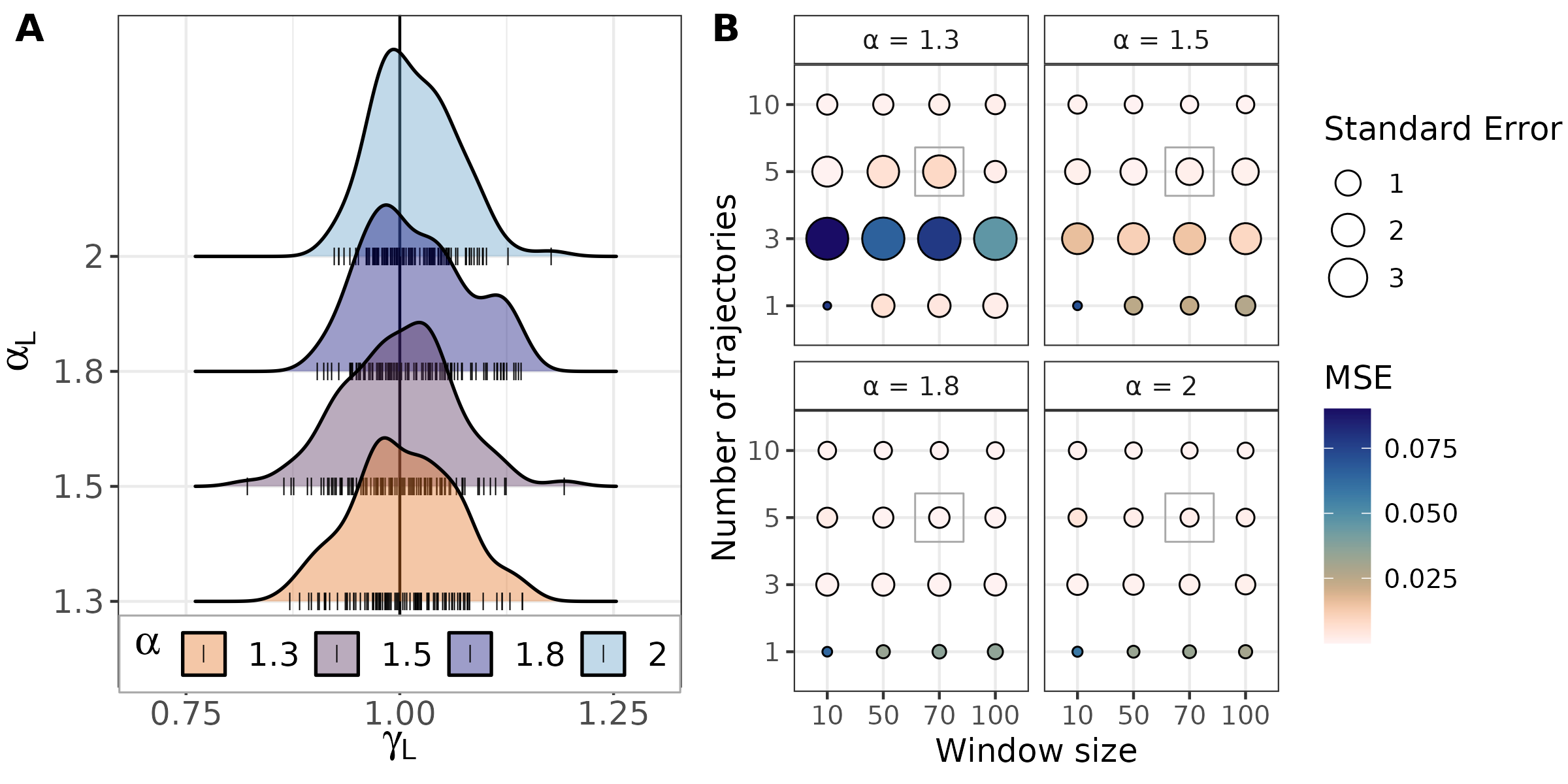}
    \caption{Performance of algorithm. \textbf{A}: Accuracy of algorithm to estimate $\gamma_N$ for different values of $\alpha_N$ (True $\gamma_N = 1$, sample size is 200. \textbf{B} Accuracy of estimating $\gamma_X$ in dependence of window size per trajectory and number of trajectories used. Multiplying both gives total amount of data points per estimation. Color represents mean square error $\frac{| \gamma_{\textnormal{true}} - \gamma_{\textnormal{est}} | ^2}{\gamma_{\textnormal{true}}}$ and size represents the standard error $\frac{\sigma_{\gamma_{\textnormal{est}}}}{\gamma_{\textnormal{true}}}$. Estimation time increases with number of data points (not shown). Grey box indicates parameters chosen for simulation.}
    \label{fig:benchmark_gamma}
\end{figure}

For the equilibrium runs we perform 100 independent estimations of $\gls{gamma}_{X}$ for each combination of $\alpha$ and $k$. As our goal here was to confirm our theoretical findings, we use 5 independent trajectories for each estimation to improve accuracy at reasonable computational costs (see Figure~\ref{fig:benchmark_gamma}). All parameters used in the simulations are also given in Table~\ref{tab:simulation}. To reduce the influence of stochasticity on our estimations, we use the same noise sequence across the range of $k$ within each estimation and the same random seed to generate noise sequences for different $\alpha$ (see Figure \ref{fig:trajectories} for an illustration of the latter). 

For the non-equilibrium runs, we simulated 15 trajectories for each value of $\alpha$. After reaching equilibrium, we varied $k$ from 5 to 0 in steps of 0.0001. We estimated $\gamma_X$ every 150 time steps, using 300 data points.

\begin{table}[]
    \centering
    \caption{Overview of simulation parameters}
    \begin{tabular}{l l r}
    \hline
    Parameter & Equilibrium & Non-Equilibrium \\
    \hline
    \multicolumn{2}{l}{\textbf{Euler-Maruyama}} \\
    Time step $\Delta t$ & 0.004 & 0.004\\
    Number of timesteps  & 10000 & 10000 + 50000 \\
    $\gamma_N$ $\sigma$ & 0.1 & 0.1\\
    $\alpha$ & 2.0, 1.8, 1.5, 1.3 & 2.0, 1.8, 1.5, 1.3 \\
    \multicolumn{2}{l}{\textbf{Estimation of $\mathbf{\gamma_X}$}}\\
    Sample size    &  100 & 15\\
    Number of data points & 70 X 5 trajectories & 300\\
    Values of $k$ & 100, 10, 1, 0.1, 0.01  &  5 to 0 by 0.0001 \\
    \hline
    \end{tabular} 
    \label{tab:simulation}
\end{table}

\subsection{Equilibrium simulations}

\begin{figure}
    \centering
    \includegraphics[width=\textwidth]{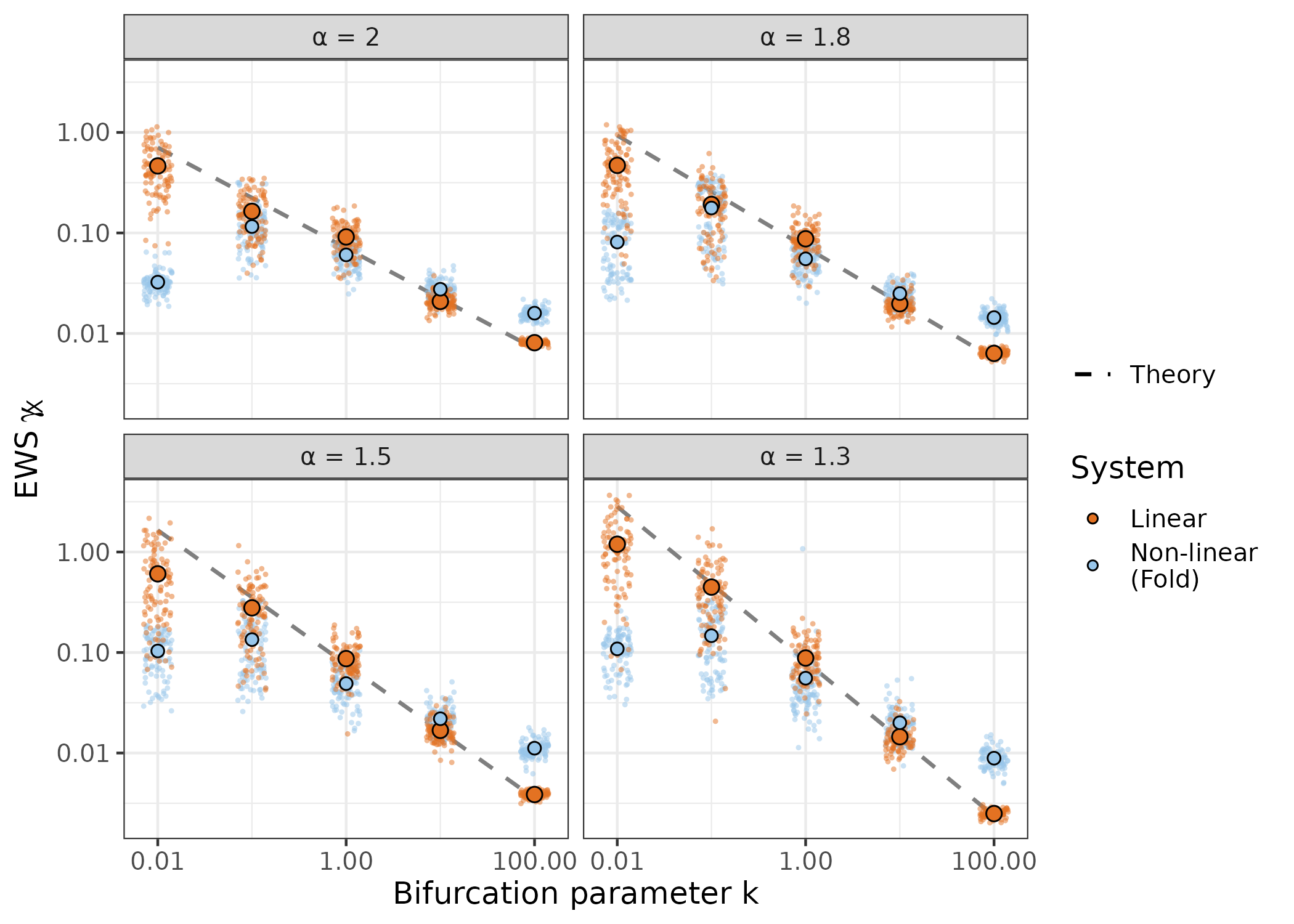}
    \caption{Equilibrium simulations: Estimation of $\gls{gamma}_{X}$ as a function of $k$ for the (linear) Ornstein-Uhlenbeck process \eqref{eq:oup} (orange) and non-linear fold \eqref{eq:fold} (blue) on a log-log scale. Small dots represent 100 individual estimations, large dots the mean value across the whole sample. Grey lines give the theoretical result \eqref{eq:gammaX}.}
    \label{fig:gamma_eq}
\end{figure}

Our simulations of the Ornstein-Uhlenbeck process confirm the theoretical relationship between $k$ and $\gamma_X$ (Figure~\ref{fig:gamma_eq}).  Accuracy is highest for large values of $k$; the smaller $k$, the higher the variability between independent estimations. However, the mean across simulations corresponds to theoretical values for all $k$ and $\alpha$, only deviating slightly very close to the bifurcation. 

In the non-linear case, we see similar patterns of increasing variability for lower values of $k$ and $\alpha$. Mean values align with theory for medium values of $k$ but not very far or very close to the bifurcation. This is expected as the linearization \eqref{eq:linearization} neglects higher-order terms, which become more important as we approach the bifurcation point. Nevertheless, we observe a strong increase in $\gamma_X$ up until $k = 0.1$, confirming the theoretical suitability of $\gamma_X$ as an early warning sign across all simulated $\alpha$ for a wide range of $k$. 
 
As expected, estimating $\gamma_X$ from trajectories produces more noisy results, with individual trajectories exhibiting large jumps in $\gamma_X$, especially for smaller $\alpha$ due to large jumps of the underlying process. 

The mean across trajectories fits the theoretical value well at the start of the simulation but begins to deviate more and more as the simulation progresses. This is consistent with theory as we are leaving the equilibrium case and the system takes longer to reach equilibrium again as we move towards a bifurcation. However, $\gamma_X$ continues to increase. An exception is the linear case for $\alpha$s of 1.3 and 1.5, where we see a stagnation or even decline of the mean trajectory very close to the bifurcation ($k < 1$). Importantly in the non-linear case, this is not the case and we observe a steady increase in $\gamma_X$ for the whole range of $k$ and all $\alpha$ in both the mean and the majority of individual trajectories. This confirms the practical suitability of $\gamma_X$ as an early warning sign of an approaching bifurcation in more application-oriented situations.

\begin{figure}
    \centering
    \includegraphics[width=\textwidth]{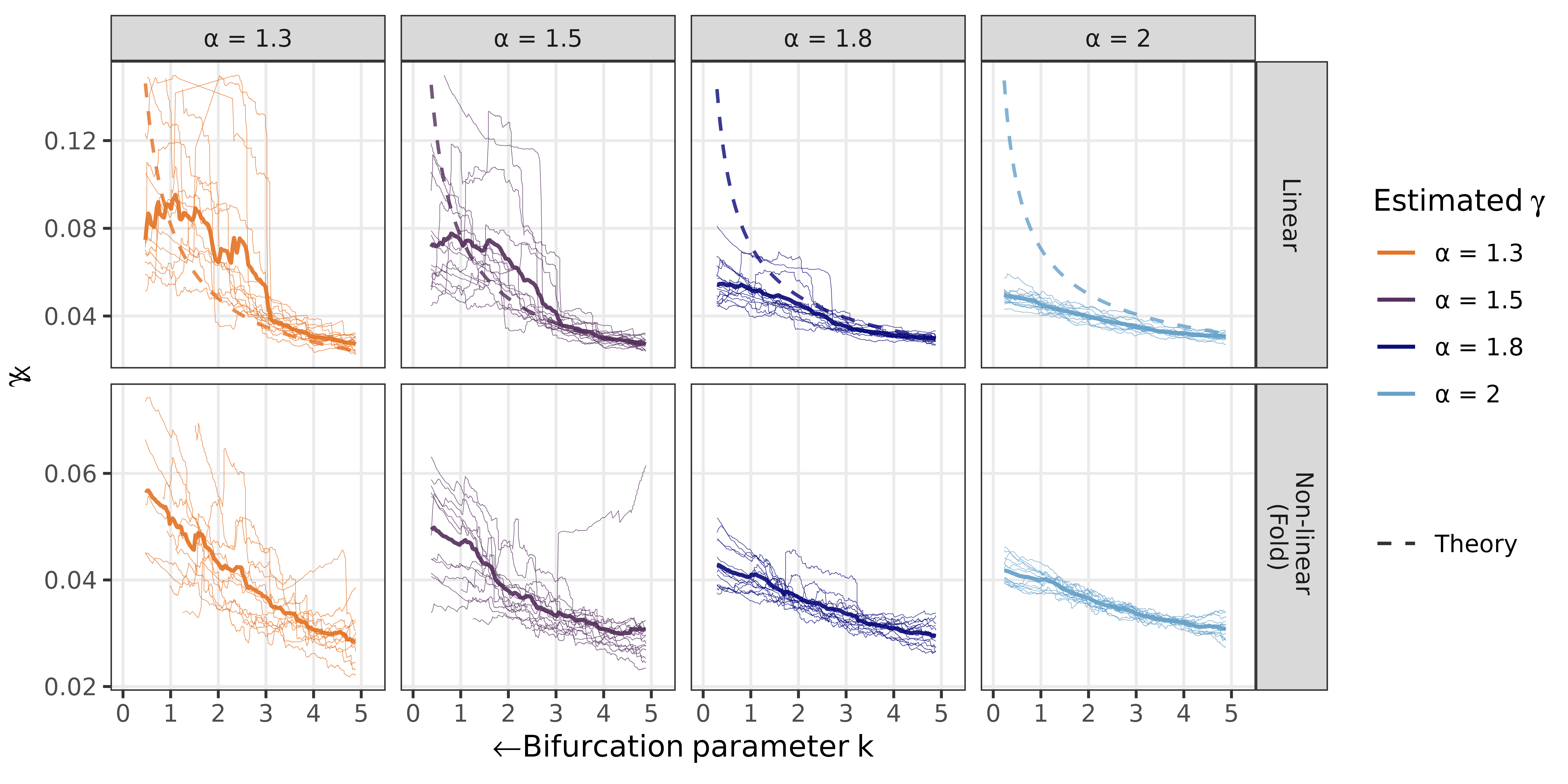}
    \caption{Non-equilibrium simulations: Estimation of $\gls{gamma}_{X}$ on transient trajectories produced by the linear and non-linear systems. Thin lines represent 15 individual estimations, thick lines the mean value across the whole sample. All trajectories are additionally smoothed using a moving window of 100 points. Dashed lines give theoretical results. $k$ is moved in the direction of the arrow.}
    \label{fig:gamma_neq}
\end{figure}

\section{Conclusion}

We have shown that for systems driven by $\alpha$-stable, non-Gaussian noise, the classical early warning sign of rising variance and autocorrelation are not supported by mathematical theory and its use poses the danger of spurious, false-positive results.  

To address this, we have introduced the scaling factor $\gamma_X$ as alternative, generalized early warning sign applicable to Gaussian and non-Gaussian $\alpha$-stable processes. We have laid out the necessary mathematical theory to show $\gamma_X$ is always defined and inversely scales with the bifurcation parameter, much in the same way as the variance does in the Gaussian case. 

Our simulations confirmed our theoretical results and showed that $\gamma_X$ can be estimated from few trajectories with sufficient accuracy. Additionally, our results generalize well to the non-linear, non-equilibrium case we would usually find in applications. 

Estimating the parameters of an $\alpha$-stable distribution is a common exercise and algorithms are readily available in relevant programming languages. While being computationally more expensive than variance estimation, it still provides an easy-to-use method that works with a limited amount of data points available. This provides good conditions for applying $\gamma_X$ to more complex and real-world data streams in the future.

With $\alpha$-stable models again gaining traction in climate and tipping point research, we thus hope our results will contribute to their further understanding and use.

\FloatBarrier

\section{Appendix}

We are interested in the statistical properties of the process $X$. We thus formulate its characteristic function

\begin{equation}
     \gls{charf}_{X}(u) = \mathbb{E}[e^{iuX(t)}] 
\end{equation}

and, using Eq. \eqref{eq:stableoupsolution} and initial conditions $X_0 = 0$ obtain

\begin{equation}
     \gls{charf}_{X}(u) = \mathbb{E}[e^{iu\int_0^t e^{-\gls{bifurcation}(t-s)}dN(s)}]
\end{equation} \label{eq:solutionincharf}

Making use of the It\^{o}-Integral

\begin{align}
    \gls{charf}_{X}(u) & =  \lim_{L\to\infty} \mathbb{E}[ e^{iu e^{-\gls{bifurcation}t} \sum_{j=1}^L e^{\gls{bifurcation}s_{j}^{(L)}} (N(s^{(L)}_{j+1}) - N(s^{(L)}_{j}))}] \\
    & = \lim_{L\to\infty} \prod_{j=1}^L \mathbb{E} [ e^{iu e^{-\gls{bifurcation}(t - s_{j}^{(L)})}(N(s^{(L)}_{j+1}) - N(s^{(L)}_{j}))} ]
\end{align}

If we map the last expression to \eqref{eq:charfL} and take $u e^{-\gls{bifurcation}(t - s_{j}^{(L)})}$ as our new Fourier parameter we obtain 

\begin{align}
    \gls{charf}_{X}(u) & =  \lim_{L\to\infty} \prod_{j} e^{(s^{(L)}_{j+1} - s^{(L)}_{j})(\gls{gamma}_{N}^{\gls{alpha}}|e^{-\gls{bifurcation}(t - s_{j}^{(L)})}|^{\gls{alpha}})} \\
    & =  \lim_{L\to\infty}  e^{\sum_{j} (s^{(L)}_{j+1} - s^{(L)}_{j})(\gls{gamma}_{N}^{\gls{alpha}}|e^{-\gls{bifurcation}(t - s_{j}^{(L)})}|^{\gls{alpha}})} \\
    & = e^{-\int_0^t \gls{gamma}_{N}^{\gls{alpha}}|e^{-\gls{bifurcation}(t - s)}|^{\gls{alpha}} dt}\\
    & = e^{- \frac{\gls{gamma}_{N}^{\gls{alpha}}}{\alpha k} (1 - e^{-\alpha k t}) |u|^{\alpha}}
\end{align}

\FloatBarrier

\printbibliography

\end{document}